\documentclass[journal]{IEEEtran}  

\IEEEoverridecommandlockouts

\usepackage{graphicx}
\usepackage{epsfig} 
\usepackage{times} 
\usepackage{amsmath}
\usepackage{amssymb}  
\usepackage{color}
\usepackage{amsfonts}
\usepackage{subfigure}
\usepackage{multirow}
\usepackage{multicol}
\usepackage{siunitx}
\usepackage{overpic}
\usepackage{mathrsfs}  

\usepackage{booktabs}
\usepackage{threeparttable}

\usepackage{epstopdf}
\newtheorem{proposition}{Proposition}
\newtheorem{theorem}{Theorem}

\usepackage{xspace}

\begin{document}

\title{\LARGE \bf PDE-Based Feedback Control of Freeway Traffic Flow \\via Time-Gap Manipulation of ACC-Equipped Vehicles}

\author{Nikolaos Bekiaris-Liberis and Argiris I. Delis

\thanks{N. Bekiaris-Liberis is with both the Departments of Electrical \& Computer Engineering and Production Engineering \& Management, Technical University of Crete, Chania, Greece, 73100. Email address: \texttt{nikos.bekiaris@dssl.tuc.gr}.}

\thanks{A. I. Delis is with the Department of Production Engineering \& Management, Technical University of Crete, Chania, Greece, 73100. Email address: \texttt{adelis@science.tuc.gr}.}}

%C. Roncoli is with the Department of Built Environment, School of Engineering, Aalto University, 02150 Espoo, Finland. Email address: \texttt{claudio.roncoli@aalto.fi}.}}

\maketitle
\pagestyle{headings}

\markboth{Submitted to IEEE Transactions on Control Systems Technology on November 20, 2018}{Bekiaris-Liberis \& Delis}

\begin{abstract}
We develop a control design for stabilization of traffic flow in congested regime, based on an Aw-Rascle-Zhang-type (ARZ-type) Partial Differential Equation (PDE) model, for traffic consisting of both ACC-equipped (Adaptive Cruise Control-equipped) and manual vehicles. The control input is the value of the time-gap setting of ACC-equipped and connected vehicles, which gives rise to a problem of control of a $2\times 2$ nonlinear system of first-order hyperbolic PDEs with in-domain actuation. The feedback law is designed in order to stabilize the linearized system, around a uniform, congested equilibrium profile. Stability of the closed-loop system under the developed control law is shown constructing a Lyapunov functional. Convective stability is also proved adopting an input-output approach. The performance improvement of the closed-loop system under the proposed strategy is illustrated in simulation, also employing four different metrics, which quantify the performance in terms of fuel consumption, total travel time, and comfort.
%We design a feedback law for stabilization of traffic flow in congested regime based on an Aw-Rascle-Zhang-type (ARZ-type) Partial Differential Equation (PDE) model, for traffic consisting of both ACC-equipped (Adaptive Cruise Control-equipped) and manual vehicles. The control input is the value of the time-gap setting of ACC-equipped and connected vehicles, which gives rise to a problem of control of a $2\times 2$ nonlinear system of first-order hyperbolic PDEs with in-domain actuation. The feedback law is designed in order to stabilize the linearized system, around a uniform, congested equilibrium profile. Stability of the closed-loop system under the developed control law is shown constructing a Lyapunov functional. Convective stability is also proved adopting an input-output approach. The performance improvement of the closed-loop system under the proposed strategy is illustrated in simulation, also employing four different metrics, which quantify the performance in terms of fuel consumption, total travel time, and comfort.
\end{abstract}
\section{Introduction}
%\subsection{Motivation}
Although traffic congestion may be unavoidable nowadays, due to the continuous increase in the number of vehicles and in the traffic demand, some of its ramifications may be alleviated employing real-time traffic control strategies \cite{varaia}. Among other reasons, certain traffic flow instability phenomena, such as, for example, stop-and-go waves, are some of the causes of traffic congestion's negative consequences on fuel consumption, total travel time, drivers' comfort, and safety \cite{kesting}. One promising avenue to traffic flow stabilization is the development of control design tools that exploit the capabilities of automated and connected vehicles \cite{diakaki}, while retaining the distributed nature of traffic flow dynamics. It is the aim of this paper to develop a feedback law for traffic flow stabilization utilizing a PDE traffic flow model and exploiting the capabilities of ACC-equipped and connected vehicles.

%\baselineskip=1.8\normalbaselineskip

%costly but althogut anodivaldble, we may alleviate some fo its negative ffects 

%stabilizstion of traffic nbeacus it has negativeneffects in fuel  bla bla

%at the same time expoitng capabilities of tecnologyc

%\subsection{Literature}
Since second-order, PDE traffic flow models (i.e., systems that incorporate two PDE states, one for traffic density and one for traffic speed) constitute realistic descriptions of the traffic dynamics, capturing important phenomena, such as, for example, stop-and-go traffic, capacity drop, etc. \cite{dellis num}, \cite{goatin2}, \cite{ngoduy}, boundary control designs are recently developed for such systems \cite{bayen}, \cite{karafyllis-nikos}, \cite{goatin2}, \cite{yu1}, \cite{yu2}, \cite{prieur1}, \cite{zhang1} some of which are based on techniques originally developed for control of systems of hyperbolic PDEs, such as, for example, \cite{vazquez2}, \cite{argomedo}, \cite{herty}, \cite{lamare}, \cite{litrico}, \cite{prieurintro}, \cite{vazquez1}. Even though simpler, first-order traffic flow models, in conservation law or Hamilton-Jacobi PDE formulation, are also important for modeling purposes. For this reason, PDE-based control design techniques exist for this class of systems as well \cite{bek-bayen}, \cite{blandin}, \cite{claudel1}, \cite{delle}, \cite{goatin1}, \cite{lihor}. 

While most of the above PDE-based traffic control techniques rely on traditional implementation means such as, ramp metering and variable speed limits, more rare are PDE-based, traffic flow control methodologies that exploit connected and automated vehicles capabilities. In particular, \cite{darbha}, \cite{horowitz} develop control designs via in-domain manipulation of acceleration of ACC-equipped vehicles, considering traffic with only automated vehicles, and \cite{piacentini}, \cite{koga} develop control designs via speed manipulation of an autonomous vehicle. Furthermore, although in microscopic simulation it is reported that it may be beneficial for traffic flow, to appropriately manipulate in real time the ACC settings of vehicles already equipped with an ACC feature \cite{kesting gap}, \cite{sackel}, \cite{spiliopoulou}, the problem of systematic feedback control design via time-gap manipulation hasn't, heretofore, been tackled from a PDE viewpoint. 

In this work, we design a feedback control strategy for stabilization of traffic flow in congested regime, manipulating the time-gap setting of vehicles equipped with ACC and utilizing a control-oriented, ARZ-type model with ACC (which is shown to possess certain important traffic flow-theoretic properties). The control strategy is developed for the linearized system around a uniform, congested equilibrium profile, which is proved to be open-loop unstable. Due to the presence (on average) of a certain penetration rate of ACC-equipped vehicles in a given freeway stretch, the traffic flow control problem is recast to the problem of stabilization of a $2\times 2$ linear system of first-order, heterodirectional hyperbolic PDEs with in-domain actuation. The closed-loop system under the proposed controller is shown to be exponentially stable (in $C^1$ norm), constructing a Lyapunov functional. We further study, employing an input-output approach, an additional stability property of the closed-loop system, namely convective stability, which is important from a traffic control point of view as it guarantees the non-amplification of speed perturbations, as these propagate upstream. The benefits in traffic flow of employing the proposed strategy are illustrated in simulation, also including the quantification of the performance improvement in terms of various indices, measuring total travel time, fuel consumption, and comfort level.

%already existing commercail vehicles   providing a practiclalt feasible 

%adanvces in communiticatin

%supression of oscillations

%here we exploit the avialability of acc equipped vehicles and aim at changin its settings, providing a practiclalt feasible 

%Utilizing this model we design a feedback law that provides the desired time-gap setting for ACC-equipped vehicles. The feedback law is shown to  

%\subsection{Organization}
We start in Section \ref{sec model} presenting a control-oriented traffic flow model for congested and mixed (i.e., consisting of both manual and ACC-equipped vehicles) traffic. In Section \ref{control design} we introduce our feedback design and in Section \ref{section stability} we prove the stability and convective stability of the closed-loop system. The effectiveness of the proposed strategy is validated in simulation in Section \ref{simul}. Concluding remarks and future research directions are provided in Section \ref{conlus}.

%The model has several traffic flow-oriented properties, which makes it suitable for control purposes, such as, for example, it is anisotropic, bla bla. 

{\em Notation:} For scalar functions $u\in L^p\left(0,D\right)$, where $p$ is a positive integer and $D>0$, we define the norm $\|u\|_{p}=\left(\int_0^Du(x)^pdx\right)^{\frac{1}{p}}$ as well as the weighted norm $\|u\|_{\mu,p}=\left(\int_0^De^{p\mu x}u(x)^pdx\right)^{\frac{1}{p}}$, for $\mu\neq0$. For $u\in C[0,D]$ we denote $\|u\|_{C}=\max_{x\in[0,D]}\left|u(x)\right|=\lim_{p\to+\infty}\|u\|_{p}$ and $\|u\|_{\mu,C}=\max_{x\in[0,D]}\left|e^{\mu x}u(x)\right|=\lim_{p\to+\infty}\|u\|_{\mu,p}$. For $u\in C^1[0,D]$ we define $\|u\|_{C^1}=\|u\|_{C}+\|u'\|_{C}$ and, respectively, we define $\|u\|_{\mu,C^1}=\|u\|_{\mu,C}+\|u'\|_{\mu,C}$. For a signal $f\in \mathscr{L}_p$ we define its temporal norm $\|f\|_{\mathscr{L}_p}=\left(\int_0^{+\infty}|f(t)|^pdt\right)^{\frac{1}{p}}$, for $p<+\infty$, and $\|f\|_{\mathscr{L}_p}=\sup_{t\geq0}|f(t)|$, for $p=+\infty$. %The Laplace transform of signal $f$, $t\geq0$, is denoted by $F(s)$.%=\mathcal{L}\left\{f(t)\right\}$. 

\section{ARZ-Type Model with ACC\\ in Congested Regime}
\label{sec model}
\subsection{Description of the model}
We consider the following system
\setlength{\arraycolsep}{0pt}\begin{eqnarray}
\rho_t(x,t)&=&-\rho_x(x,t)v(x,t)-\rho(x,t)v_x(x,t)\label{sys1}\\
v_t(x,t)&=&-\left(v(x,t)+\rho(x,t)\frac{\partial V_{\rm mix}\left(\rho(x,t),h_{\rm acc}(x,t)\right)}{\partial \rho}\right)\nonumber\\
&&\times v_x(x,t)+\frac{V_{\rm mix}\left(\rho(x,t),h_{\rm acc}(x,t)\right)-v(x,t)}{\tau_{\rm mix}}\label{speed}\\
{q_{\rm in}}&=&\rho(0,t)v(0,t)\label{boundary qin}\\
v_t(D,t)&=&\frac{V_{\rm mix}\left(\rho(D,t),h_{\rm acc}(D,t)\right)-v(D,t)}{\tau_{\rm mix}},\label{sysn}
\end{eqnarray}\setlength{\arraycolsep}{5pt}where
\setlength{\arraycolsep}{0pt}\begin{eqnarray}
V_{\rm mix}\left(\rho,h_{\rm acc}\right)&=&\tau_{\rm mix}\left(\frac{\alpha}{\tau_{\rm acc}}V_{\rm acc}\left(\rho,h_{\rm acc}\right)+\frac{1-\alpha}{\tau_{\rm m}}V_{\rm m}\left(\rho\right)\!\right)\label{vmix}\\
V_{\rm acc}\left(\rho,h_{\rm acc}\right)&=&\frac{1}{h_{\rm acc}}\left(\frac{1}{\rho}-L\right),\quad\rho_{\rm min}<\rho<\frac{1}{L}\label{vacc}\\
V_{\rm m}\left(\rho\right)&=&\frac{1}{h_{\rm m}}\left(\frac{1}{\rho}-L\right),\quad\rho_{\rm min}<\rho<\frac{1}{L}\label{vm}\\
\tau_{\rm mix}&=&\frac{1}{\frac{\alpha}{\tau_{\rm acc}}+\frac{1-\alpha}{\tau_{\rm m}}},\label{tmix}
\end{eqnarray}\setlength{\arraycolsep}{5pt}$\rho$ is traffic density, $0< v\leq v_{\rm f}$ is traffic speed, with $v_f$ being some maximum achievable speed (or, free-flow speed), $D>0$ is length of a given freeway stretch, $L>0$ is average effective length of each vehicle, $\alpha\in[0,1]$ is the percentage of ACC-equipped vehicles with respect to total vehicles, $\rho_{\rm min}>0$ is the lowest value for density for which the model is accurate\footnote{One may view $\rho_{\rm min}$ as the critical density that corresponds to a minimum possible time-gap (see Section \ref{subtra}).}, $t\geq 0$ is time, $x\in\left[0,D\right]$ is spatial variable, $q_{\rm in}>0$ is a constant external inflow, $\tau_{\rm acc}$, $\tau_{\rm m}$ $>0$ are the time constants of the ACC-equipped and manual vehicles, respectively, $h_{\rm m}>0$ is the time-gap of manual vehicles, and $h_{\rm acc}>0$ is the time-gap of ACC-equipped vehicles, which is the control input.

%\footnote{Fundamental diagrams of this form correspond to the so-called constant time-gap (CTG) policy, see, e.g., \cite{bose2}, \cite{darbha}, \cite{horowitz}.}, respectively.

\subsection{Traffic flow-oriented properties of the model}
\label{subtra}
\begin{figure}[h]
\centering
\includegraphics[width=0.85\linewidth]{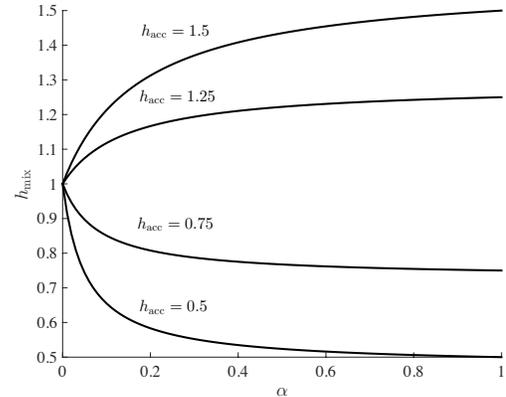}
\caption{Mixed time-gap (\ref{eff head}) for $h_{\rm m}=1$, $\frac{ \tau_{\rm acc}}{\tau_{\rm m}}=0.1$, and four different values of ${h}_{\rm acc}$, as a function of the penetration rate $\alpha$.}
\label{fig1}
\end{figure}

The motivation for model (\ref{sys1})--(\ref{sysn}) is the following. First, note that equation (\ref{sys1}) is the conservation-of-vehicles equation, as $q=\rho v$ is the traffic flow. Equation (\ref{speed}) is the speed equation, which is inspired by the speed dynamics of the ARZ model \cite{zhang}. In fact, the ARZ model may be viewed as both a model of traffic flow dynamics for traffic with only manual vehicles \cite{zhang} as well as a model for traffic flow dynamics with only ACC-equipped vehicles \cite{darbha} (Section 3.2). In particular, for fixed time-gaps of ACC-equipped (and manual) vehicles, when $\alpha=1$ (only ACC-equipped vehicles exist) or $\alpha=0$ (only manual vehicles exist) the model reduces to the ARZ model with fundamental diagram given by (\ref{vacc}) or (\ref{vm}) (which corresponds to the so-called constant time-gap policy, see, e.g., \cite{bose2}, \cite{darbha}, \cite{horowitz}), respectively. However, to account for the case of mixed traffic, i.e., when both manual and ACC-equipped vehicles are present, we define a new equilibrium (fundamental diagram) relation for speed as in (\ref{vmix}), which is also written as 
\begin{eqnarray}
V_{\rm mix}\left(\rho,h_{\rm acc}\right)=\frac{1}{h_{\rm mix}\left(h_{\rm acc}\right)}\left(\frac{1}{\rho}-L\right), \label{alt vmix}
\end{eqnarray}
where the effective (or, mixed) time-gap is defined as
\begin{eqnarray}
h_{\rm mix}\left(h_{\rm acc}\right)=\frac{\alpha+\left(1-\alpha\right)\frac{ \tau_{\rm acc}}{\tau_{\rm m}}}{\alpha+\left(1-\alpha\right)\frac{ \tau_{\rm acc}}{\tau_{\rm m}}\frac{{h}_{\rm acc}}{h_{\rm m}}}{h}_{\rm acc}.\label{eff head}
\end{eqnarray}We show in Fig. \ref{fig1} the mixed time-gap as a function of the penetration rate $\alpha$ for $\frac{ \tau_{\rm acc}}{\tau_{\rm m}}=0.1$, ${h}_{\rm m}=1$, and four different values for the time-gap of ACC-equipped vehicles. 

\begin{figure}[h]
\centering
\includegraphics[width=0.85\linewidth]{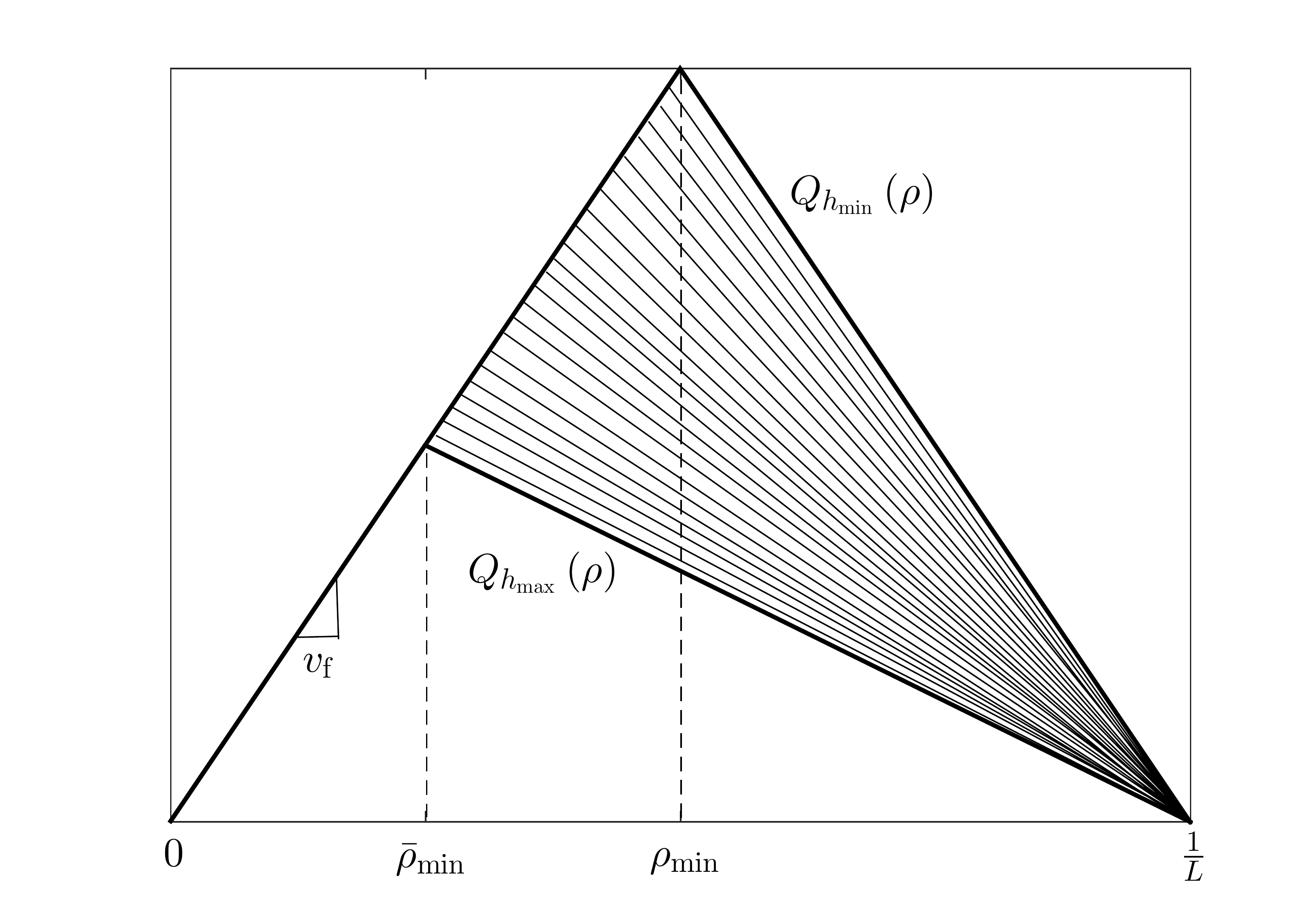}
\caption{Different fundamental diagrams (\ref{vacc}) for  $h_{\rm acc}\in\left[h_{\rm min},h_{\rm max}\right]$.}
\label{fig FD}
\end{figure}

In the present work, we restrict our attention to congested regime, and thus, it is sufficient to define only the right part (i.e., for $\frac{1}{L}>\rho>\rho_{\min}$) of fundamental diagrams (\ref{vacc}), (\ref{vm}). However, one may utilize any appropriate extension for the left part (i.e., for $0<\rho\leq\rho_{\rm min}$), such as, for example, a fundamental diagram that corresponds to a constant (free-flow) speed. We show in Fig. \ref{fig FD} an example of potentially meaningful fundamental diagrams (\ref{vacc}) for different (but fixed) values of $h_{\rm acc}$. 
Specifically, $Q_{h_{\rm min}}$ is fundamental diagram that corresponds to some minimum possible time-gap, say $h_{\rm min}$, which is related to $\rho_{\rm min}$ via $h_{\rm min}=\frac{\frac{1}{\rho_{\rm min}}-L}{v_{\rm f}}$, defined as\footnote{Although $Q_{h_{\rm min}}$ is not differentiable at $\rho_{\rm min}$, one could obtain a differentiable approximation of the original fundamental diagram by adding an $\epsilon$-layer around the critical density and defining $Q_{h_{\rm min}}$ for $\rho\in\left[\rho_{\rm min}-\epsilon,\rho_{\rm min}+\epsilon\right]$ properly. Since we don't deal with free-flow conditions, in order to not distract the reader with additional technical details we don't discuss this further.}
\begin{eqnarray}
Q_{h_{\rm min}}\left(\rho\right)=\left\{\begin{array}{ll}v_{\rm f}\rho,&0\leq\rho\leq\rho_{\rm min}\\\frac{1}{h_{\rm min}}\left(1-L\rho\right),&\rho_{\rm min}<\rho\leq\frac{1}{L}\end{array}\right..
\end{eqnarray}
The fundamental diagram that corresponds to some maximum possible time-gap, say $h_{\rm max}$, is defined respectively as $Q_{h_{\rm max}}\left(\rho\right)=v_{\rm f}\rho$ for $0\leq\rho\leq\bar{\rho}_{\rm min}$ and $Q_{h_{\rm max}}\left(\rho\right)=\frac{1}{h_{\rm max}}\left(1-L\rho\right)$ for $\bar{\rho}_{\rm min}<\rho\leq\frac{1}{L}$, where $h_{\rm max}$ satisfies $h_{\rm max}=\frac{\frac{1}{\bar{\rho}_{\rm min}}-L}{v_{\rm f}}$ (for further details on realistic values of $h_{\rm max}$ and $h_{\rm min}$ that may appear in practice see Section \ref{simul} as well as, e.g., \cite{nok}, \cite{spiliopoulou}). Every other fundamental diagram (\ref{vacc}) that may appear, for different time-gaps of ACC-equipped vehicles within the interval $\left[h_{\rm min},h_{\rm max}\right]$, lies between $Q_{h_{\rm max}}$ and $Q_{h_{\rm min}}$. Furthermore, for a given penetration rate, since relation (\ref{eff head}) implies that $\min\left\{h_{\rm acc},h_{\rm m}\right\}\leq h_{\rm mix}\leq\max\left\{h_{\rm acc},h_{\rm m}\right\}$, for all $\alpha\in[0,1]$, it follows that whenever $\min\left\{h_{\rm acc},h_{\rm m}\right\}\geq h_{\rm min}$ and $\max\left\{h_{\rm acc},h_{\rm m}\right\}\leq h_{\rm max}$, all of the possible mixed fundamental diagram relations (\ref{alt vmix}) that may appear, for any $\alpha\in[0,1]$, lie between $Q_{h_{\rm max}}$ and $Q_{h_{\rm min}}$, and hence, as long as $\rho>\rho_{\rm min}=\frac{1}{L+v_{\rm f}h_{\rm min}}$ the mixed fundamental diagram relation (\ref{alt vmix})  corresponds to congested traffic. 

In addition, from (\ref{tmix}) it follows that $\min\left\{\tau_{\rm acc},\tau_{\rm m}\right\}\leq \tau_{\rm mix}\leq\max\left\{\tau_{\rm acc},\tau_{\rm m}\right\}$, for all $\alpha\in[0,1]$, and hence, when $\tau_{\rm acc}<\tau_{\rm m}$, which is typically the case in practice, $\tau_{\rm mix}$ is a decreasing function of $\alpha$. Since for given values of $v_{\rm f}$ (dependent, for example, on the specific freeway stretch) and $L$, the requirements $\min\left\{h_{\rm acc},h_{\rm m}\right\}\geq h_{\rm min}$ and $\max\left\{h_{\rm acc},h_{\rm m}\right\}\leq h_{\rm max}$ guarantee that $0<V_{\rm mix}\left(\rho,h_{\rm acc}\right)< v_{\rm f}$, for all $\alpha\in[0,1]$ and  $\rho_{\rm min}<\rho<\frac{1}{L}$, we obtain a speed equation that may serve as a reasonable model for speed in case of mixed traffic in congested conditions.

Since we are concerned with the case of congested traffic conditions we restrict our attention in a nonempty, connected open subset $\Omega$ of the set 
\noindent
$\bar{\Omega}=\left\{\left(v,\rho,h_{\rm acc}\right)\in\mathbb{R}^3:0< v< v_{\rm f},\frac{1}{L+v_{\rm f}h_{\rm min}}<\rho<\frac{1}{L},\right.$ $\left.h_{\rm min}\leq h_{\rm acc}\leq h_{\rm max}\right\}$, such that $v\!+\!\rho\frac{\partial V_{\rm mix}\left(\rho,h_{\rm acc}\right)}{\partial \rho}\!<\!0$, for all $\alpha\in[0,1]$, whenever $\left(v,\rho,h_{\rm acc}\right)\in{\Omega}$\footnote{Provided that $\max\left\{h_{\rm acc},h_{\rm m}\right\}\!\leq\! h_{\rm max}$, this holds true, for instance, for all $\left(v,\rho,h_{\rm acc}\right)\in\bar{\Omega}$ in the special (but quite restrictive) case where $v_{\rm f}\!\!\leq\!\!\frac{L}{h_{\rm max}}$.}, see, e.g., \cite{bayen}, \cite{kesting}. In fact, from (\ref{alt vmix}), it is evident that (\ref{sys1})--(\ref{sysn}) is a $2\times 2$ system of first-order hyperbolic PDEs with (real and distinct) eigenvalues given by $\lambda_1=v$ and $\lambda_2=v+\rho\frac{\partial V_{\rm mix}\left(\rho,h_{\rm acc}\right)}{\partial \rho}=v-\frac{1}{h_{\rm mix}\left(h_{\rm acc}\right)\rho}$, which implies that information propagates forward with traffic flow at the traffic speed, whereas speed information travels backward. Thus, model (\ref{sys1})--(\ref{sysn}) is anisotropic, see, e.g., \cite{zhang}.

\subsection{Boundary conditions of the system}
The boundary condition (\ref{boundary qin}) at the inlet of the considered freeway stretch implies that the flow at the entrance of the freeway stretch is equal to some external inflow with value $q_{\rm in}$. To obtain a realistic downstream boundary condition as well as to obtain a well-posed system we impose the dynamic boundary condition (\ref{sysn}), which implies free downstream traffic conditions, see also, e.g., \cite{karafyllis-nikos}. This is reasonable, even under congested conditions (consider, for example, the case where at the outlet of the considered stretch there is the end of a tunnel or the end of high-curvature or the end of an upgrade, etc.).

%see also, e.g.

%To obtain free downstream traffic conditions, which is a reasonable thing to ask even under congested conditions (consider, for example, the case where at the outlet of the considered highway stretch there is the end of a tunnel or the end of high-curvature or the end of an upgrade, etc.), as well as to obtain a well-posed system we impose the dynamic boundary condition (\ref{sysn}), see also, e.g., \cite{karafyllis-nikos}.

\subsection{Equilibria of the system}
The equilibria of system (\ref{sys1})--(\ref{sysn}) dictated by a constant inflow $q_{\rm in}$ as well as a constant, steady-state time-gap for ACC-equipped vehicles, say $\bar{h}_{\rm acc}$, which results in a steady-state mixed time-gap given by
\begin{eqnarray}
\bar{h}_{\rm mix}=\frac{\alpha+\left(1-\alpha\right)\frac{ \tau_{\rm acc}}{\tau_{\rm m}}}{\alpha+\left(1-\alpha\right)\frac{ \tau_{\rm acc}}{\tau_{\rm m}}\frac{\bar{h}_{\rm acc}}{h_{\rm m}}}\bar{h}_{\rm acc},\label{gap steady}
\end{eqnarray}
are uniform and satisfy
\begin{eqnarray}
\bar{v}&=&\frac{q_{\rm in}}{\bar{\rho}},\label{rhoeq}
\end{eqnarray}
as well as the fundamental diagram relation 
\begin{eqnarray}
\frac{1}{\bar{\rho}}-L&=&\bar{h}_{\rm mix}\bar{v}.\label{veq}
\end{eqnarray}
To see this, first note that relations (\ref{sys1}) and (\ref{boundary qin}) imply that the equilibrium values for $\rho$ and $v$, say $\rho^{\rm e}$ and $v^{\rm e}$, respectively, satisfy $\rho^{\rm e}(x)v^{\rm e}(x)=q_{\rm in}$, for all $x\in[0,D]$. From (\ref{speed}) and (\ref{sysn}) it then follows, using (\ref{alt vmix}), that the equilibrium profile of the speed satisfies the following ODE in $x$
\begin{eqnarray}
{v^{\rm e}}'(x)&=&-\frac{1}{\tau_{\rm mix}}\frac{v^{\rm e}(x)+\frac{L}{\bar{h}_{\rm mix}-\frac{1}{q_{\rm in}}}}{v^{\rm e}(x)},\label{eq ode}
\end{eqnarray}
with final condition $v^{\rm e}\left(D\right)=-\frac{L}{\bar{h}_{\rm mix}-\frac{1}{q_{\rm in}}}$. Thus, 
\begin{eqnarray}
v^{\rm e}(x)=\frac{L}{\frac{1}{q_{\rm in}}-\bar{h}_{\rm mix}}=\bar{v}, \quad \mbox{for all $x\in[0,D]$},\label{eqvd}
\end{eqnarray}
 which can be seen noting that $v^{\rm e}=\bar{v}$ is an equilibrium of (\ref{eq ode}). In order to guarantee that $\rho_{\rm min}<\bar{\rho}<\frac{1}{L}$, $\forall \alpha\in[0,1]$, which also implies from (\ref{rhoeq}), (\ref{veq}) that $0<\bar{v}<\frac{1}{\bar{h}_{\rm mix}}\left(\frac{1}{\rho_{\rm min}}-L\right)\leq v_{\rm f}$, we require that time-gaps and inflow are such that relation $0<q_{\rm in}<\frac{v_{\rm f}h_{\rm min}}{h_{\rm max}\left(L+v_{\rm f}h_{\rm min}\right)}$ holds.

\section{Control Design for the Linearized System}
\label{control design}
\subsection{Linearization and diagonalization of the system}
We start defining the error variables $\tilde{\rho}(x,t)=\rho(x,t)-\bar{\rho}$, $\tilde{v}(x,t)=v(x,t)-\bar{v}$, and $\tilde{h}_{\rm acc}(x,t)=h_{\rm acc}(x,t)-\bar{h}_{\rm acc}$. Linearizing system (\ref{sys1})--(\ref{sysn}) around the uniform, congested equilibrium profile we get
\setlength{\arraycolsep}{2pt}\begin{eqnarray}
\tilde{\rho}_t(x,t)+\bar{v}\tilde{\rho}_x(x,t)+\bar{\rho}\tilde{v}_x(x,t)&=&0\label{nonre1}\\
\tilde{v}_t(x,t)-c_4\tilde{v}_x(x,t)&=&-c_1\tilde{\rho}(x,t)-c_2\tilde{v}(x,t)\nonumber\\
&&-c_3\tilde{h}_{\rm acc}(x,t)\\
\tilde{\rho}(0,t)+c_5\tilde{v}(0,t)&=&0\\
\tilde{v}_t(D,t)&=&-c_1\tilde{\rho}(D,t)-c_2\tilde{v}(D,t)\nonumber\\
&&-c_3\tilde{h}_{\rm acc}(D,t),\label{nonren}
\end{eqnarray}\setlength{\arraycolsep}{5pt}where $c_1=\frac{1}{\bar{\rho}^2\tau_{\rm mix}\bar{h}_{\rm mix}}$, $c_2=\frac{1}{\tau_{\rm mix}}$, $c_3=\frac{\alpha}{\tau_{\rm acc}\bar{h}^2_{\rm acc}}\left(\frac{1}{\bar{\rho}}-L\right)$, $c_4=\frac{L}{\bar{h}_{\rm mix}}$, and $c_5=\frac{\bar{\rho}}{\bar{v}}$. Defining $\tilde{z}(x)=e^{\frac{c_2x}{\bar{v}}}\left(\tilde{\rho}(x)+\bar{h}_{\rm mix}\bar{\rho}^2\tilde{v}(x)\right)$ and noting that $c_2-c_1\bar{h}_{\rm mix}\bar{\rho}^2=0$, we re-write system (\ref{nonre1})--(\ref{nonren}) in diagonal form as
\setlength{\arraycolsep}{0pt}\begin{eqnarray}
\!\!\!\!\!\tilde{z}_t(x,t)+\bar{v}\tilde{z}_x(x,t)&=&-e^{\frac{c_2x}{\bar{v}}}\bar{h}_{\rm mix}\bar{\rho}^2c_3\tilde{h}_{\rm acc}(x,t)\label{reim1}\\
\!\!\!\!\!\tilde{v}_t(x,t)-c_4\tilde{v}_x(x,t)&=&-c_1e^{-\frac{c_2x}{\bar{v}}}\tilde{z}(x,t)-c_3\tilde{h}_{\rm acc}(x,t)\label{bad}\\
\!\!\!\!\!\tilde{z}(0,t)&=&-L\frac{\bar{\rho}^2}{\bar{v}}\tilde{v}(0,t)\\
\!\!\!\!\!\tilde{v}_t(D,t)&=&-c_1e^{-\frac{c_2}{\bar{v}}D}\tilde{z}(D,t)-c_3\tilde{h}_{\rm acc}(D,t).\label{reimn}
\end{eqnarray}\setlength{\arraycolsep}{5pt}%, which follows from (\ref{c1def}), (\ref{defc2}).

%\begin{eqnarray}
%\tilde{w}_t(x,t)+\bar{v}\tilde{w}_x(x,t)&=&-\bar{h}_{\rm mix}\bar{\rho}^2c_1\tilde{w}(x,t)-\bar{h}_{\rm mix}\bar{\rho}^2\left(c_2-c_1\bar{h}_{\rm mix}\bar{\rho}^2\right)\tilde{v}(x,t)-\bar{h}_{\rm mix}\bar{\rho}^2c_3\tilde{h}_{\rm acc}(x,t)\\
%\tilde{v}_t(x,t)-c_4\tilde{v}_x(x,t)&=&-c_1\tilde{w}(x,t)-\left(c_2-c_1\bar{h}_{\rm mix}\bar{\rho}^2\right)\tilde{v}(x,t)-c_3\tilde{h}_{\rm acc}(x,t)\\
%\tilde{w}(0,t)&=&-L\frac{\bar{\rho}^2}{\bar{v}}\tilde{v}(0,t)\\
%\tilde{v}_t(D,t)&=&-c_1\tilde{w}(D,t)-\left(c_2-c_1\bar{h}_{\rm mix}\bar{\rho}^2\right)\tilde{v}(D,t)-c_3\tilde{h}_{\rm acc}(D,t).
%\end{eqnarray}

% defining $\tilde{z}(x,t)=\tilde{w}(x,t)e^{\frac{c_2x}{\bar{v}}}$ we get

%\subsection{Motivation for feedback control}

\subsection{Control law}
In addition to improving performance, feedback control is needed because system (\ref{reim1})--(\ref{reimn}) for ${h}_{\rm acc}=\bar{h}_{\rm acc}$ is unstable, as it is shown in the next proposition whose proof can be found in Appendix A.
\begin{proposition}
\label{prop1}
System (\ref{reim1})--(\ref{reimn}) is exponentially unstable in open-loop.
\end{proposition}

%MAKES SENSE FROM TRAFFIC POINT OF VIEW? READ WHAT KARAFYLLIS SAID TO MIROSLAV

%SAY SOMEHWEHRE THAT IS ACC INFOR AONLY NEEDED BECAUSE OF THE SPATIONAL DERIVATIV WRT TO SPEED AS RAJAGOPAL SAY 

%read what karafyllis SAID ABOUT SPOURIOS EQULIBIRA THAT DO NOT MAKE SENSE THE NONUNIFORM ONES

%\begin{proof}
%\end{proof}

The control law is chosen as
\setlength{\arraycolsep}{0pt}\begin{eqnarray}
{h}_{\rm acc}(x,t)&=&\bar{h}_{\rm acc}+\frac{1}{c_3}\left(-c_1e^{-\frac{c_2x}{\bar{v}}}\tilde{z}(x,t)+k\tilde{v}(x,t)\right)\\
&=&\bar{h}_{\rm acc}+\frac{1}{c_3}\left(-c_1\tilde{\rho}(x,t)+\left(k-c_2\right)\tilde{v}(x,t)\right),\label{controller}
\end{eqnarray}\setlength{\arraycolsep}{5pt}with $k>0$ being arbitrary, which gives
\setlength{\arraycolsep}{0pt}\begin{eqnarray}
\tilde{z}_t(x,t)+\bar{v}\tilde{z}_x(x,t)&=&c_2\tilde{z}(x,t)-ke^{\frac{c_2x}{\bar{v}}}\bar{h}_{\rm mix}\bar{\rho}^2\tilde{v}(x,t)\label{tilz}\\
\tilde{v}_t(x,t)-c_4\tilde{v}_x(x,t)&=&-k\tilde{v}(x,t)\label{speed1}\\
\tilde{z}(0,t)&=&-L\frac{\bar{\rho}^2}{\bar{v}}\tilde{v}(0,t)\\
\tilde{v}_t(D,t)&=&-k\tilde{v}(D,t).\label{speed2}
\end{eqnarray}\setlength{\arraycolsep}{5pt}From the closed-loop system (\ref{tilz})--(\ref{speed2}) it is evident that the feedback law aims at eliminating the source term in (\ref{bad}), which may cause instability due to a feedback connection between the states $\tilde{z}$ and $\tilde{v}$, while rendering the $\tilde{v}\left(D\right)$ subsystem exponentially stable (and autonomous).

%SAY ABOUT MEASURMENT REQUIREMENTS AND REFER TO MY ESTIMATION PAPERS 

%MENTION CONNECTED VEHICLES ETC

%say hOW IT IS IMPLENENTED IN PRACGICE WITH CONNECTED VEHICLES FROM DIAMANTI

Taking into account that the traffic system operates in congested regime, the operating point of the controller, as this is seen via the steady-state time-gap for ACC-equipped vehicles $\bar{h}_{\rm acc}$, may vary considering, for example, safety, comfort, or total travel time criteria. For instance, in cases in which safety is a primary goal, the time-gap $\bar{h}_{\rm acc}$ may take large values (which implies that $\bar{h}_{\rm mix}$ also takes large values, according to (\ref{gap steady})), whereas, when comfort is of significant importance, then no action (e.g., as recommendation to drivers of ACC-equipped vehicles or as direct manipulation of the ACC settings of individual vehicles) may be taken (in order to not disrupt the driver) from the controller for imposing the value of $\bar{h}_{\rm acc}$, which implies that the driver alone may set the value for the time-gap $\bar{h}_{\rm acc}$, see, e.g., \cite{spiliopoulou}. Moreover, it may be beneficial, from a total travel time point of view, the time-gap $\bar{h}_{\rm acc}$ to take large values, since, for given inflow, lower steady-state densities may be achieved (via the achievement of higher steady-state speeds), as it can be seen from relations (\ref{rhoeq}), (\ref{eqvd}). We consider a specific scenario and further discuss about the choice of $\bar{h}_{\rm acc}$ (as well as of $h_{\rm m}$) in Section~\ref{simul}.

%\footnote{Note that this implies that $\bar{h}_{\rm mix}$ also takes large values since $\bar{h}_{\rm mix}$ is an increasing function of $\bar{h}_{\rm acc}$, for given penetration rate $\alpha$ and time-gap of manually driven vehicles $h_{\rm m}$, as it is evident from relation (\ref{gap steady}).}

%MAYBE JUST SAY TOTAL TRAVEL TIME THUS SMALLER DENSITY NOT SPEED

In practice, under a vehicle-to-infrastructure (V2I) communication paradigm, the control authority may implement the proposed strategy either as time-gap recommendations to drivers of ACC-equipped vehicles or via direct manipulation of the ACC settings of such vehicles, see, e.g., \cite{spiliopoulou}. Furthermore, the developed feedback law, given in the simple formulae (\ref{controller}), requires measurements of the average speed and density (or, equivalently, average speed and flow, via the flow definition $q=\rho v$, in case flow measurements are available instead) throughout the spatial domain. This information could be obtained by the central control authority via utilization of connected vehicles\footnote{Besides ACC-equipped vehicles, a connected vehicle may be any vehicle able to exchange information with the central monitoring and control unit.} reports (e.g., reporting speed, position, or other information) as well as measurements from fixed detectors and, potentially, also employing certain traffic state estimation methodologies, see, e.g., \cite{bek est1}, \cite{claudel2}, \cite{bek est2}, \cite{wang-est}. 

%see alos section simulaton and papers 

%Since we are in congested regime and since it is evident from (\ref{gap steady}) that $\bar{h}_{\rm mix}$ is an increasing function of $\bar{h}_{\rm acc}$, we choose the minimum possible steady-state time-gap to increase the discharge flow \cite{spiliopoulou} as well as, namely 
%\begin{eqnarray}
%\bar{h}_{\rm acc}= \frac{\frac{1}{\rho_{\rm min}}-L}{v_{\rm f}}.\label{con steady}
%\end{eqnarray}

%WHAT THE DRIVER OF ACC VEHICLES CHOOSES IN CONGESTION??? I ASK BECAUSE IN MY SIMULATIONS I WILL CONSIDER COMFORT SO NO RECOMMENDATION SO THE VALUE FOR ACC IS WHAT IS SETS USUALLY BY HUMAN DRIVERS.

%MAYBE SAY THAT THE STAEYD STATE DEPENDS ON THE PEENTRATION RATE AND THE STEADY STATE SETTINGS IN CONGESTION FOR MNAULAUL AND ACC (AS WELL AS THE TAUS). FOR EXMAPLE, MANUAL VEHICLES IN CONGESTION USUALLY OPERATE AT THE MINIMUM WHILE ACC VEHICLES TYPICALLY HAVE LARGER SETTINGS FOR SAFETY. ACCORIND TO, WE COULD EQUAL THE TWO AS A STRATEGY COULD EMPLOY THE MINIMUM OF THE TWO ASS SETTING FOR ACC VEHICLES WHICH IS ON AVERAGE THE MNAUL VEHICLES ONE.

\section{Stability and Convective Stability Analyses}
\label{section stability}
\subsection{Stability in $C^1$ norm}
We establish next stability in the stronger $C^1$ norm in order to guarantee additional stability properties for the closed-loop system that may be desirable from a traffic flow control viewpoint, see, e.g., \cite{horowitz}. Stability results in other norms, such as, for example, the $L^2$ norm, may be also obtained. The proofs of such results follow from the proof of the following theorem, which is provided in Appendix B. 

\begin{theorem}
\label{theorem1}
Consider a closed-loop system consisting of system (\ref{nonre1})--(\ref{nonren}) and control law (\ref{controller}). For all initial conditions $\left(\tilde{\rho}\left(\cdot,0\right),\tilde{v}\left(\cdot,0\right)\right)\in C^1[0,D]\times C^1[0,D]$, which satisfy first-order compatibility with boundary conditions, there exists a positive constant $\mu$ such that the following holds\footnote{The assumptions of Theorem \ref{theorem1} imply that $\left(\tilde{z}\left(\cdot,0\right),\tilde{v}\left(\cdot,0\right)\right)\in C^1[0,D]\times C^1[0,D]$ satisfy first-order compatibility, and thus, system (\ref{tilz})--(\ref{speed2}) exhibits a unique, classical solution such that $\tilde{z},\tilde{v}\in C^1\left([0,D]\times [0,+\infty)\right)$ (and hence, so does $\tilde{\rho}$), see, e.g., \cite{coron book}, \cite{courant}, \cite{russel}.} for all $t\geq0$
\begin{eqnarray}
\|\tilde{\rho}(t)\|_{C^1}+\|\tilde{v}(t)\|_{C^1}\leq\mu\left(\|\tilde{\rho}(0)\|_{C^1}+\|\tilde{v}(0)\|_{C^1}\right)e^{-\frac{k}{2} t}.\label{thm1est}
\end{eqnarray}
\end{theorem}

\subsection{Convective stability in $\mathscr{L}_p$, $p\in\left[1,+\infty\right]$, norm}

%i can name it convective 
%then i can say 
%i can say it is related to the string stability in micro models 
%this is why sometimes people say string convective 

%WHEN A DISTURBANCE APPEARS AT A SPECIFIC POINT

%IF THERYRE IS AN ON-RAMO AT A SPECIFIC POINT IN THE SPATIAL DOMAIN ACTING AS SINGULAR SOURCE

We study next the convective stability properties of the closed-loop system, see, e.g., \cite{bayen}, \cite{kesting}, which is a notion related to string stability of a finite platoon of vehicles, see, e.g., \cite{orosz trc}, \cite{swaroop0}. In a nutshell, convective stability in the present case guarantees that the magnitude (in $\mathscr{L}_p$ sense) of the deviation of speed as well as of its gradient, at some location (e.g., due to the presence of an unmodeled on-ramp at this specific location, acting as singular source), from the equilibrium point, decreases as the perturbation propagates backward in the spatial domain. Adopting an input-output approach we establish the following result whose proof can be found in Appendix C.

\begin{theorem}
\label{theorem2}
System (\ref{speed1}), (\ref{speed2}) is $\mathscr{L}_p$, $p\in\left[1,\infty\right]$, convectively stable in the sense that for any $0\leq x_2<x_1\leq D$ such that $\tilde{v}\left(x_1\right)\in{\mathscr{L}_p}$ and $\tilde{v}_x\left(x_1\right)\in{\mathscr{L}_p}$, the following hold
\begin{eqnarray}
%\int_0^{+\infty}\tilde{v}\left(x_2,t\right)^2dt\leq \int_0^{+\infty}\tilde{v}\left(x_1,t\right)^2dt.
\|\tilde{v}\left(x_2\right)\|_{\mathscr{L}_p}&<&\|\tilde{v}\left(x_1\right)\|_{\mathscr{L}_p}\label{est st}\\
\|\tilde{v}_x\left(x_2\right)\|_{\mathscr{L}_p}&<&\|\tilde{v}_x\left(x_1\right)\|_{\mathscr{L}_p}.\label{est st1}
\end{eqnarray}
\end{theorem}

%The explicit solution to (\ref{speed1}), (\ref{speed2}) is
%\begin{eqnarray}
%\tilde{v}\left(x,t\right)=\left\{\begin{array}{ll}\tilde{v}_0\left(c_4 t+{x}\right)e^{-kt},&t+\frac{x-D}{c_4}\leq 0\\ \tilde{v}_0\left(D\right)e^{-kt},&t+\frac{x-D}{c_4}>0\end{array}\right..
%\end{eqnarray}
%and
%\begin{eqnarray}
%\tilde{v}_x\left(x,t\right)=\left\{\begin{array}{ll}\tilde{v}_0'\left(c_4 t+{x}\right)e^{-kt},&t+\frac{x-D}{c_4}\leq 0\\ 0,&t+\frac{x-D}{c_4}>0\end{array}\right..
%\end{eqnarray}

%where it is assume ignoring the effect of the values within the interval $\frac{x_2-x_1}{c_4}\leq s\leq0$

\section{Simulation Results}
\label{simul}
\subsection{Model parameters and numerical implementation}
The parameters of system (\ref{sys1})--(\ref{sysn}) utilized in the simulation investigations are shown in Table \ref{table1}.
 \begin{table}[h]
\caption{Parameters of system (\ref{sys1})--(\ref{sysn}).}
\begin{center}
%\resizebox{1.5in}{!} {
\begin{tabular}{cc||cc||cc}   
\hline\hline
%Diffusion coefficient ($\epsilon$)&$.1$\\
$q_{\rm in}$&$1200\left(\frac{\textrm{veh}}{\textrm{h}}\right)$  &$\tau_{\rm acc}$&$2$ (s)&$h_{\rm m}$&$1$ (s)\\[2mm]
$\rho_{\rm min}$&$37\left(\frac{\textrm{veh}}{\textrm{km}}\right)$&$\alpha$&$0.15$&$D$&$1000$ (m) \\[2mm]
$\tau_{\rm m}$&$60$ (s)&$L$&$5$ (m) \\[2mm]
\hline
\end{tabular}
\label{table1}
\end{center}
\end{table}The chosen parameters are considered reasonable for a traffic flow model, see, e.g., \cite{bayen}, \cite{dellis num}, \cite{seibold}, \cite{ngoduy}. In particular, we choose a value for the time-gap of manually driven vehicles $h_{\rm m}$ that is close to reported average values of about $1.2$ s, see, e.g., \cite{nok}, \cite{spiliopoulou}, but slightly lower than this to reflect evidence that drivers may follow a preceding vehicle at smaller time-gaps in congested traffic, compared to the case of light traffic conditions, see, e.g., \cite{nok}. %Note that, at steady-state, these values satisfy all the requirements of Section \ref{sec model}.

%we choose a penetration rate of ACC-equipped vehicles of $20$\%, which is for the near future. 

For the numerical solution of the hyperbolic system (\ref{sys1})--(\ref{sysn}) in open-loop as well as under (\ref{controller}), a modified Rusanov scheme, which is an explicit finite-volume scheme of centered type with added numerical diffusion, with time and spatial discretization steps of $0.1$ s and $10$ m, respectively, is employed, see, e.g., \cite{Do}, \cite{rus}. The ODE (\ref{sysn}) that corresponds to the downstream boundary condition for the speed is numerically solved utilizing a forward Euler method with the same time step. The upstream and downstream boundary values for density and speed, respectively, are obtained from the boundary conditions (\ref{boundary qin}) and (\ref{sysn}), whereas for obtaining the ``missing" upstream and downstream boundary values for speed and density, respectively, we use fictitious cells, extrapolating the corresponding values from the interior of the domain.% whereas the upstream and downstream values for density and speed, respectively, are obtained from the boundary conditions (\ref{boundary qin}) and (\ref{sysn}).

\subsection{Controller's parameters and performance evaluation}
The operating point of the traffic system, as it is dictated by the steady-state value of the mixed time-gap according to (\ref{gap steady}), it is selected such that $\bar{h}_{\rm acc}=1.5$ s. Such a value reflects the fact that the equilibrium of the time-gap for ACC-equipped vehicles may be dictated from drivers' choices rather than from interventions of the control authority, for a control strategy that aims at minimizing controller's interventions, which may be disrupting for the driver. Consequently, we choose a value for $\bar{h}_{\rm acc}$ that is close to what drivers of ACC-equipped vehicles set in congested conditions, which is evidenced to be larger compared to manual driving in heavy traffic and which is reported to be around the selected value, see, e.g., \cite{nok}.

%aiming at retaining a good level of driver's comfort via.

%minimizing the volume of controller's interventions

%which may result, for example, from a control strategy that targets at driver's comfort, minimizing potentially distracting controller's interventions. 

 %considering driver's comfort, 
 
 %to not sitract the driver
 
 %would be the case for a steady state that reuslts from drivers choisces
 
% the equilibrium value of the time-gap for ACC-equipped vehicles may be the result of drivers' choices rather than of intervention of the control authority. 
%may set the value accodirnidn to the choice of the time-gap in congested conditions of drivers of acc equipped vehilses who tedn tend to choose a larger value for the time-gap

%The reason for this choice is that 

 %which is explained as follows. In congested conditions, 

%We consider a scenario in which driver's comfort is the primary goal of the steady-state operation of the control law. Since the control law may be implemented as recommendations to driver, in order to not distract the driver, the control authority doesn't provide anyn recommendation for the steady-state time-gap, but rather, operates with a tme-gap that the driver may select. Since in congested congestions the values are  usually larger, see, e.g., \cite{nok}, than we consider a scenario in which the statey state time gap is . 

%MANUAL AROUND SMALLER THAN USUAL  AROUND 1 AND AROUND 1.5 FOR ACC WITHOUTH RECOMMENDINGS THEM 

\begin{figure}[h]
\centering
\includegraphics[width=0.85\linewidth]{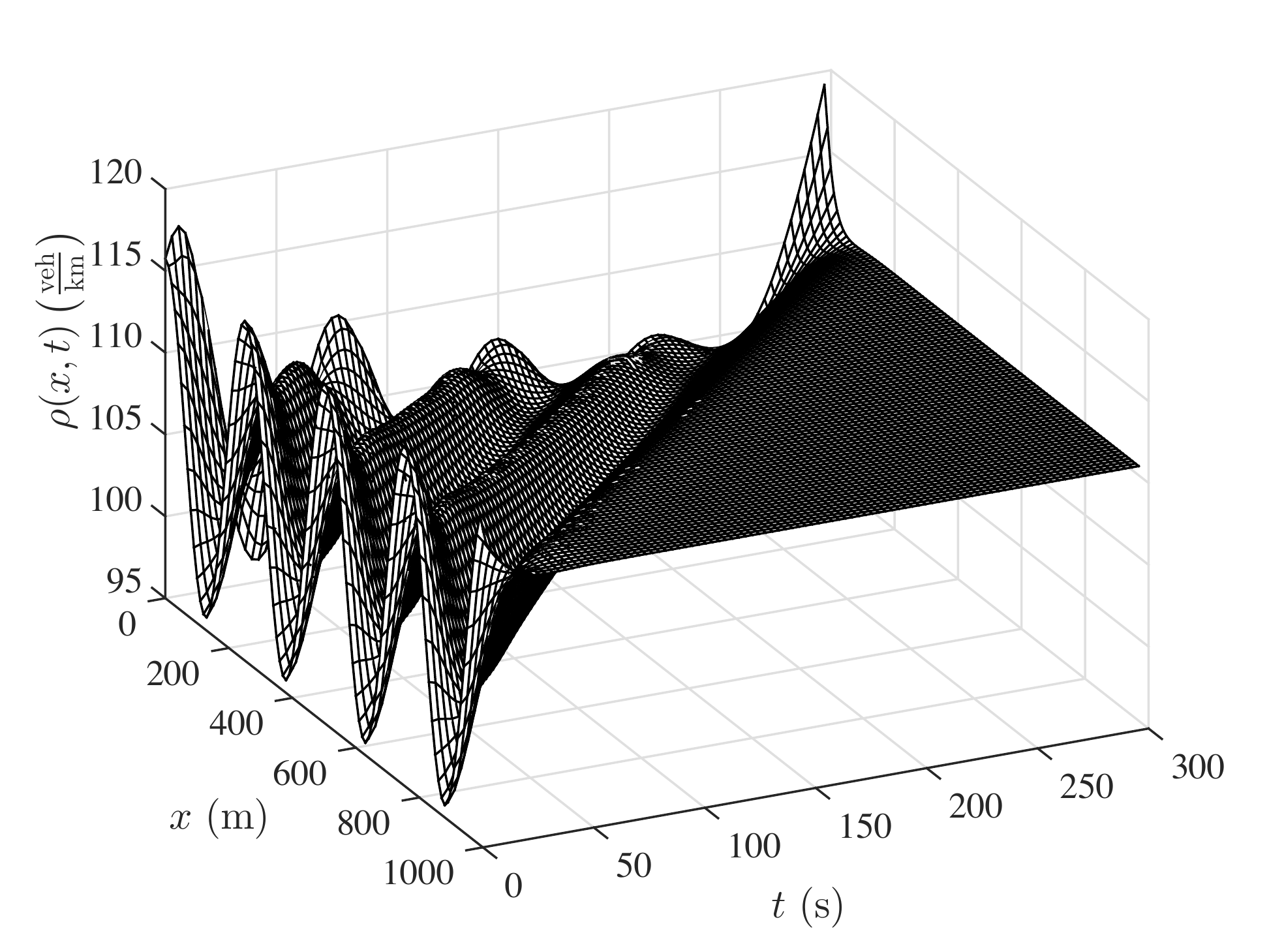}
\includegraphics[width=0.85\linewidth]{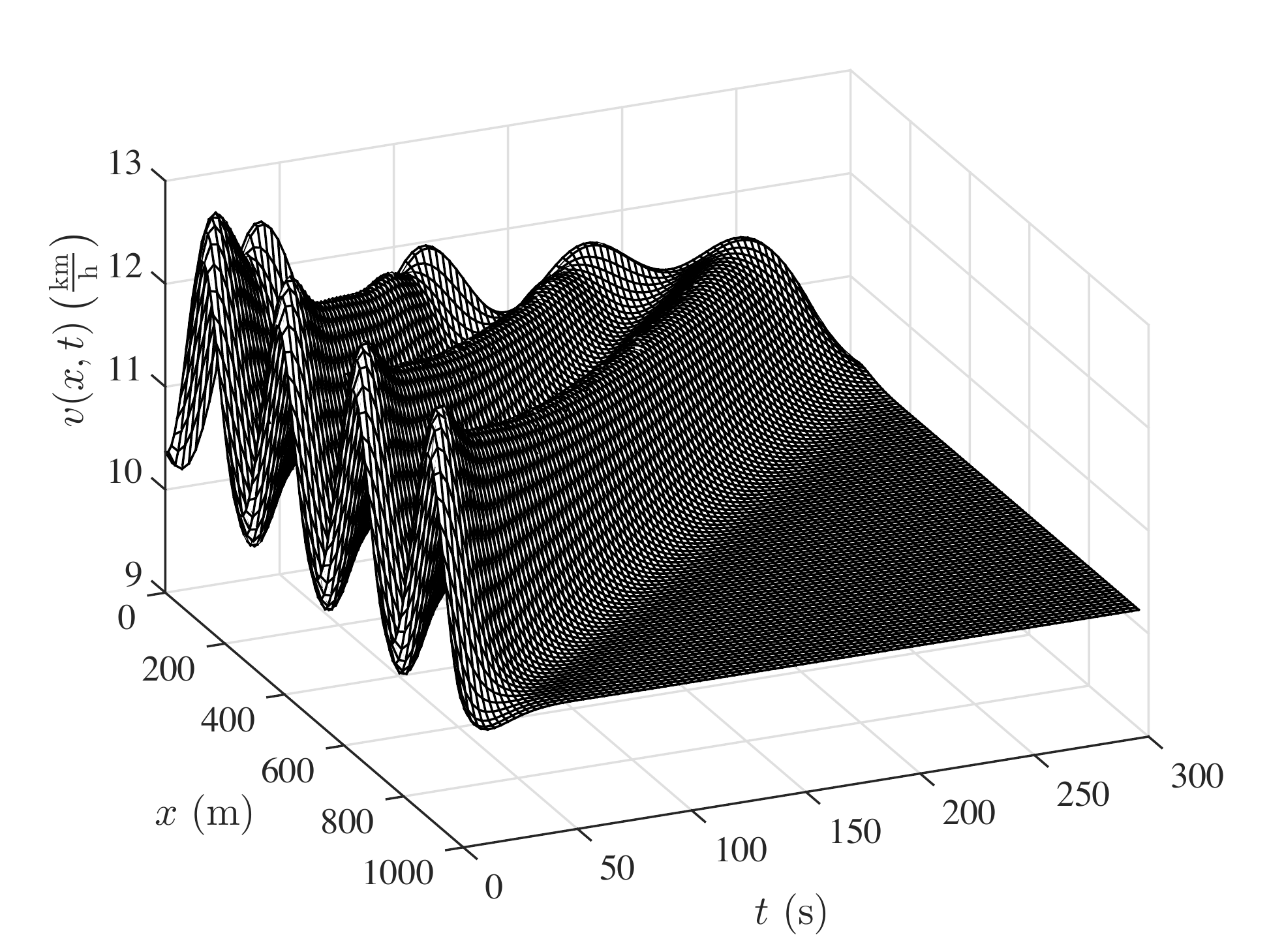}
\caption{Open-loop response of system (\ref{sys1})--(\ref{sysn}) with parameters as in Table~\ref{table1} for $\bar{h}_{\rm acc}\!\!=\!\!1.5$ and initial conditions $\rho(x,0)\!=\!\bar{\rho}+10\cos\left(\frac{8\pi x}{D}\right)$, $v(x,0)\!=\!\frac{q_{\rm in}}{\rho(x,0)}$.}
\label{fig manual}
\end{figure}

\begin{figure}[h]
\centering
\includegraphics[width=0.85\linewidth]{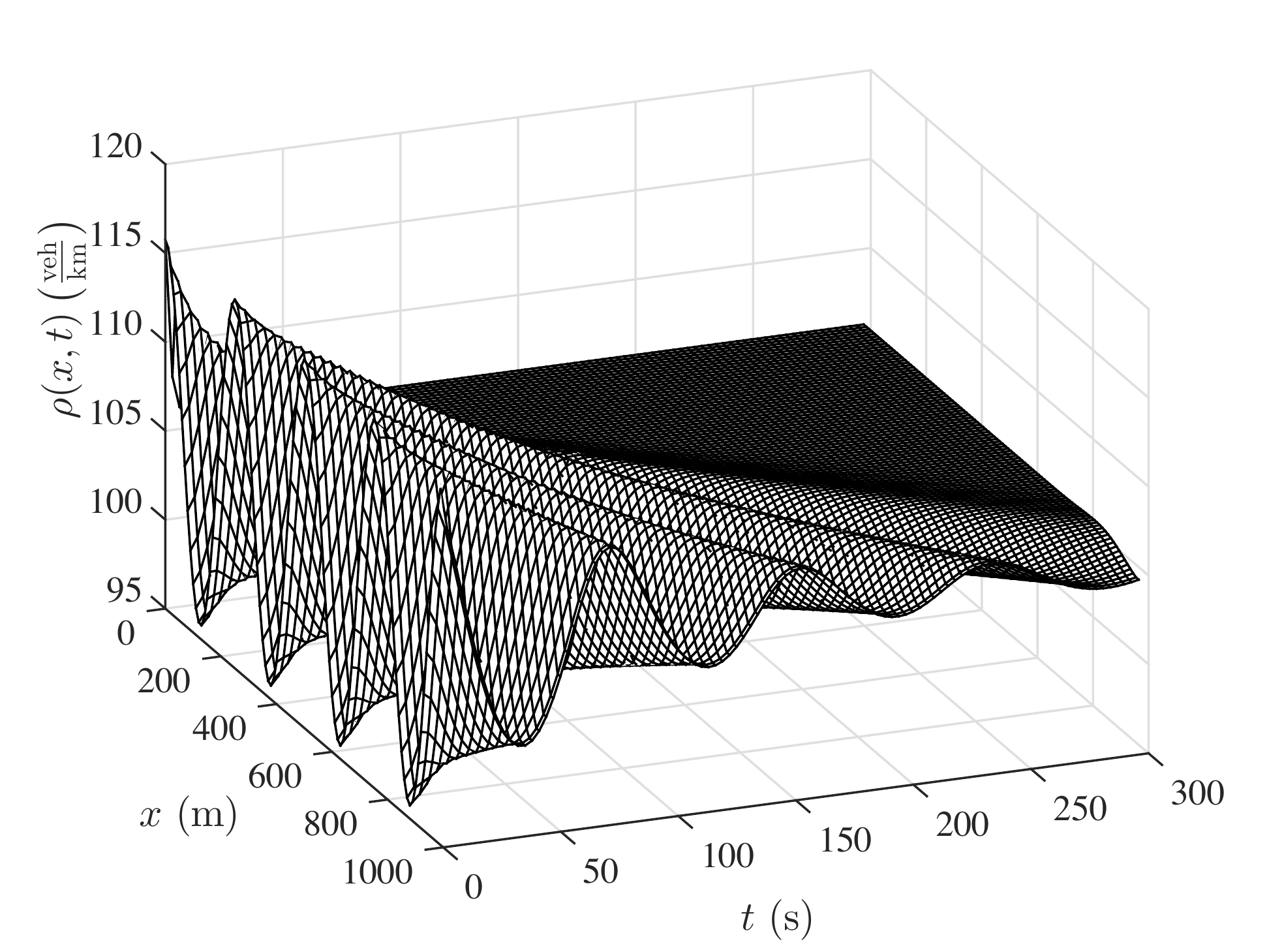}
\includegraphics[width=0.85\linewidth]{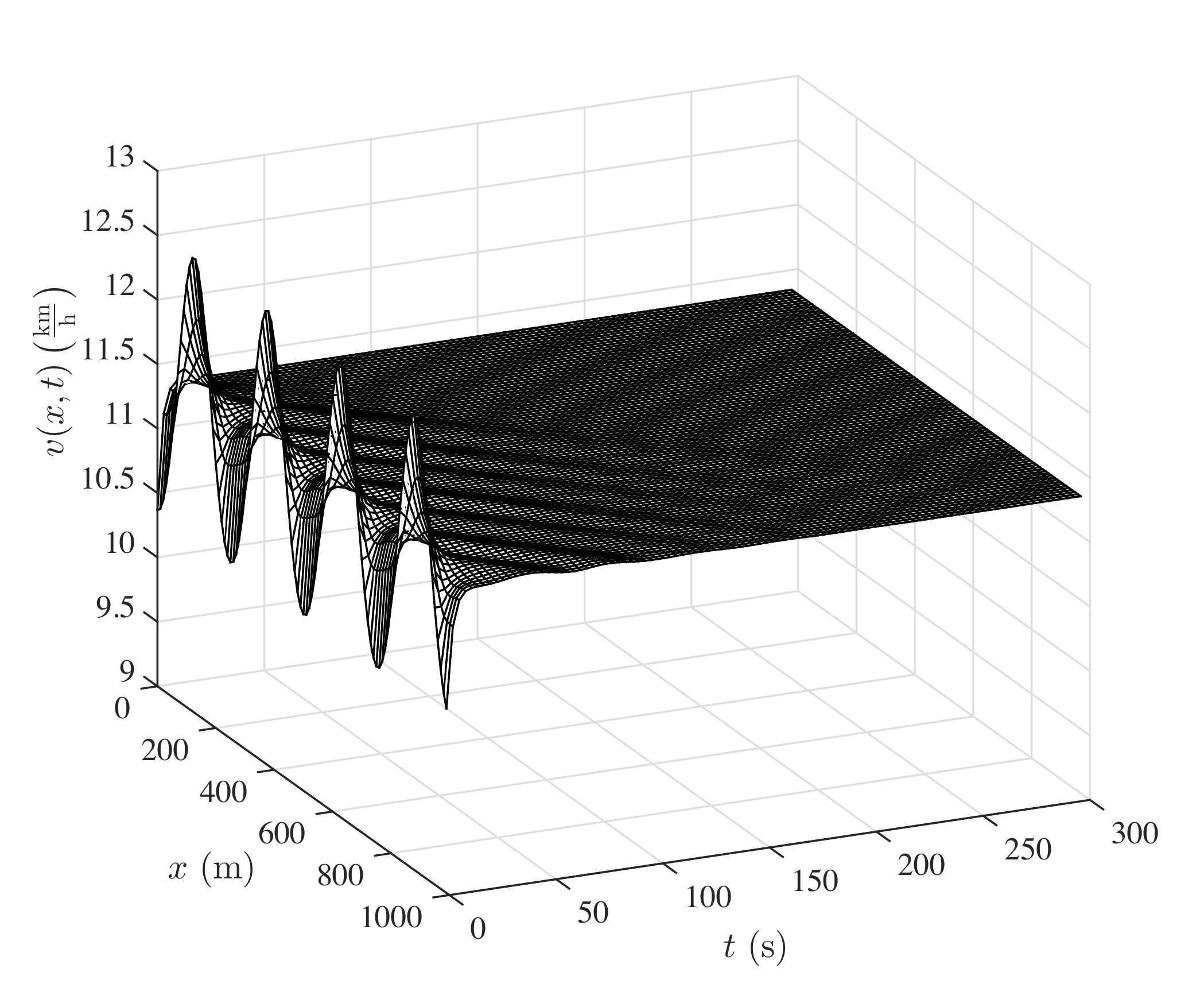}
\caption{Closed-loop response of system (\ref{sys1})--(\ref{sysn}) with parameters shown in Table \ref{table1}, under the feedback law (\ref{controller}) with $k=0.25$, for $\bar{h}_{\rm acc}=1.5$ and initial conditions  $\rho(x,0)=\bar{\rho}+10\cos\left(\frac{8\pi x}{D}\right)$, $v(x,0)=\frac{q_{\rm in}}{\rho(x,0)}$.}
\label{fig acc}
\end{figure}

\begin{figure}[h]
\centering
\includegraphics[width=0.85\linewidth]{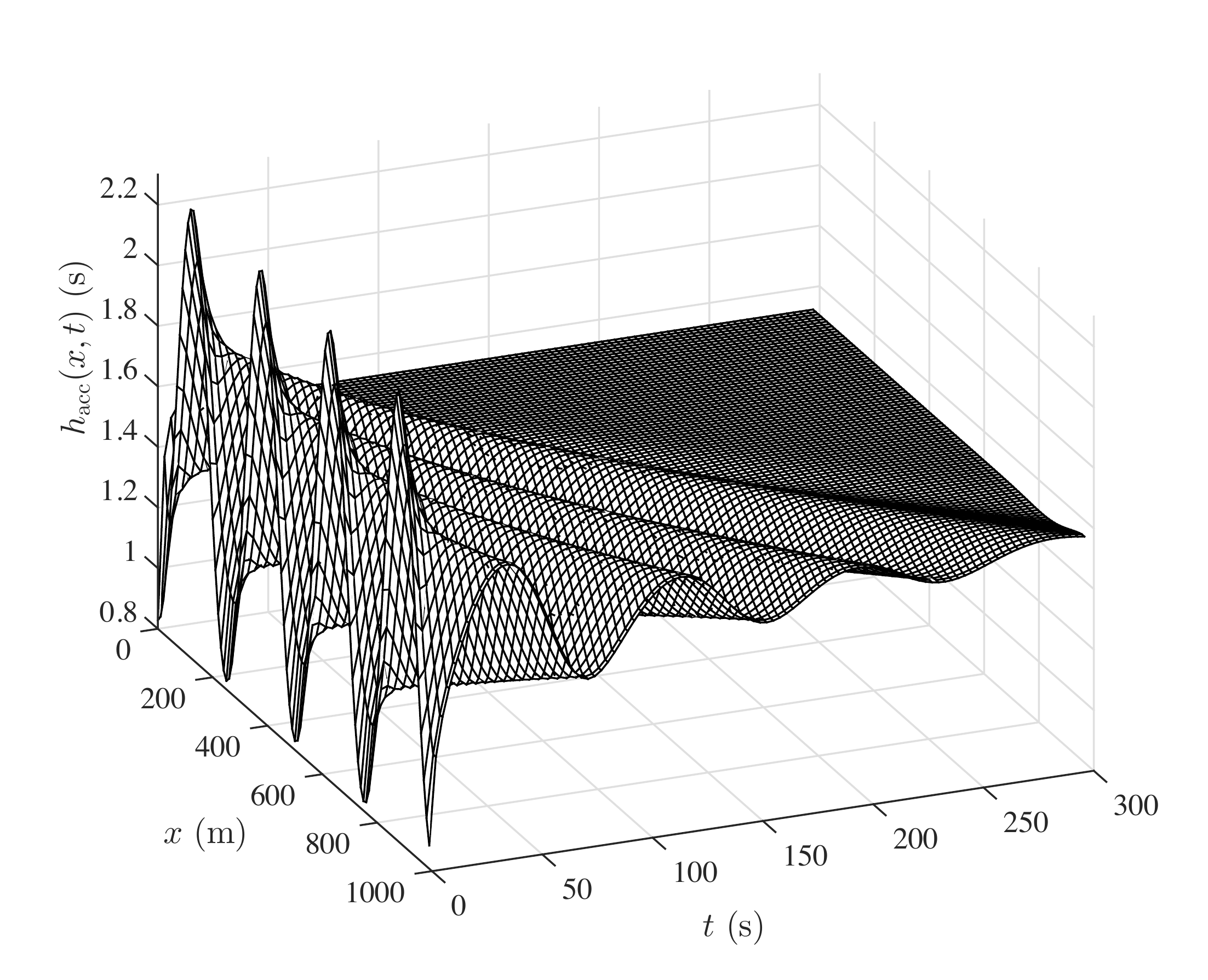}
\caption{Feedback control law (\ref{controller}) with $k=0.25$ and $\bar{h}_{\rm acc}=1.5$.}
\label{fig acc con}
\end{figure}

The steady-state values for density and speed are derived from (\ref{rhoeq}), (\ref{veq}) as $\bar{\rho}\!\!=\!\!105.8$ $\frac{\textrm{veh}}{\textrm{km}}$, $\!\bar{v}\!\!=\!\!11.35$ $\frac{\textrm{km}}{\textrm{h}}$.  We show in Fig.~\ref{fig manual} the open-loop response for initial conditions $\rho(x,0)\!=\!\bar{\rho}+\!10\cos\left(\frac{8\pi x}{D}\right)$, $v(x,0)\!\!=\!\!\frac{q_{\rm in}}{\rho(x,0)}$. From Fig. \ref{fig manual} it is evident that the open-loop response exhibits an unstable and quite oscillatory behavior. In contrast, as it is shown in Fig. \ref{fig acc}, the traffic flow is stabilized and, in particular, the oscillations (stop-and-go waves) in the speed response are considerably suppressed when the feedback law (\ref{controller}) is applied. The control effort (\ref{controller}) for $k=0.25$ $\frac{1}{\rm{s}}$ is shown in Fig. \ref{fig acc con}, from which one can also observe that the resulting values for the time-gap of ACC-equipped vehicles lie within the bounds typically implemented in ACC-equipped vehicles settings, namely, approximately within the interval $[0.8,2.2]$ s, see, e.g., \cite{nok}, \cite{spiliopoulou}.

%is oscillatory and, eventually, becomes unstable.

 %resulting eventually to an unstable response. 

%The gain $k$ of the feedback law is set to $0.25$ $\left(\frac{1}{\rm{s}}\right)$, which results in values for the time-gap of ACC-equipped vehicles within the typical bounds of about $[0.8,2.2]$, see, e.g., .

%PBSERVE THAT CHANGES THE DIRECTION OF DENSITY PRORPAGATION IN CLOSED AND OPEN

%CHECK DIRECTION OF PROPAGATION IN PLOTS

To quantify the benefits of controller (\ref{controller}) we compare the performances of the closed- and open-loop systems in terms of fuel consumption, comfort, and total travel time (TTT). We use the following performance indices, see, e.g., \cite{kesting} (Chapter~21)
\setlength{\arraycolsep}{0pt}\begin{eqnarray}
J_{{\rm fuel},1}&=&\int_0^T\int_0^D\bar{J}_{{\rm fuel},1}\left(v(x,t),a(x,t)\right)\rho(x,t)dxdt\label{cost fuel}\\
J_{{\rm comfort}}&=&\int_0^T\int_0^D\!\left(a(x,t)^2+a_t(x,t)^2\right)\rho(x,t)dxdt\label{cost comfort}\\
J_{{\rm TTT}}&=&\int_0^T\int_0^D\rho(x,t)dxdt,\label{cost ttt}
\end{eqnarray}\setlength{\arraycolsep}{5pt}where $\bar{J}_{{\rm fuel},1}\left(v,a\right)=\max\left\{0,b_0+b_1v+b_3v^3+b_4va\right\}$, $a=v_t+vv_x$, $T=350$ s, and $b_0$, $b_1$, $b_3$, $b_4$ are provided in \cite{kesting} (page 485). We further use an alternative metric for fuel consumption, namely $J_{{\rm fuel},2}$, defined in \cite{fuel delft} (Section~2.3) with parameters from \cite{ahn}, which is based on discretized, macroscopic traffic flow models. Application of the controller results in better performance in all of the metrics, as it is shown in Table~\ref{table3}. In particular, the improvement in fuel consumption and comfort is attributed to the fast homogenization of the speed field. The improvement in fuel consumption and TTT may be best appreciated taking also into account the cost of congestion, see, e.g., \cite{cost}, and the considered setup's scale. \begin{table}[h]
\caption{Performance indices (\ref{cost fuel}), (\ref{cost comfort}), (\ref{cost ttt}), and $J_{{\rm fuel},2}$ from \cite{fuel delft}.}
\begin{center}
%\resizebox{3.4in}{!} {
\begin{tabular}{|c|c|c}
\hline
Performance index&Percentage (\%) improvement with (\ref{controller})\\
\hline\hline
$J_{{\rm fuel},1}$&$3.9$\\\hline$J_{{\rm fuel},2}$&$3.3$\\\hline$J_{{\rm comfort}}$&$90$\\\hline$J_{{\rm TTT}}$&$4$\\\hline
\end{tabular}%}
\label{table3}
\end{center}
\end{table}

%IS THIS IMPROVEMENT OK?

%IS THIS SPEED TOO LOW?

%WOULD BE GOOD TO TAKE ALUE SMALLER THAN THE MANUAL GAP at some time instances

%SHOW L2?

%To choose the steady-state gap (see natasa), e.g., decreasing function of density? This may depend on $q_{\rm in}$. It may be of the form
%\begin{eqnarray}
%\bar{h}_{\rm ACC}=\left\{\begin{array}{ll}h_{\rm max},&q_{\rm in}<\\a_1-b_1q_{\rm in},&<q_{\rm in}<\\h_{\rm min},&q_{\rm in}>\end{array}\right..
%\end{eqnarray}

%is it open-loop stable? hard to tell this is a system with distriubted state delays

%when is it uncontrollable? When penetration rate is zero or when spacing is zero or density is maximum.

%disturbance rejection and 

%\section{Real-Data Validation} 

%\section{Existence and Uniqueness}

\section{Conclusions and Future Work}
\label{conlus}
We presented a control design methodology for stabilization of traffic flow in congested regime exploiting the capabilities of vehicles with ACC features and utilizing an ARZ-type model for mixed traffic. The closed-loop system, under the developed control law, was shown to be exponentially stable as well as convectively stable. The numerical investigation showed that the performance of the considered traffic system, under the proposed controller, is improved and the improvement, in terms of fuel consumption, travel time, and comfort, was quantified utilizing various performance indices.

The control design approach presented is based on a linear version of the original nonlinear system. As next step, it would be interesting to consider the nonlinear, feedback control design problem as well as to perform the analysis in a nonlinear setting considering the nonlinear closed-loop system, similarly to, e.g., \cite{coron book}, \cite{vazquez2}, \cite{hale}, \cite{karafyllis-nikos}. It would be also interesting to consider problems on circular spatial domains, as it is the case, for example, in \cite{stern}, employing a microscopic framework. %Finally, one could consider problems where there is an input delay to model various communication and actuation delays present in ACC systems, see, e.g., \cite{}.  

\setcounter{equation}{0}
\renewcommand{\theequation}{A.\arabic{equation}}
\appendices
\section*{Appendix A: Proof of Proposition \ref{prop1}}
We start computing the characteristic equation of system (\ref{reim1})--(\ref{reimn}) when $\tilde{h}_{\rm acc}\equiv 0$. One may proceed either utilizing the Laplace transform and capitalizing on the relation of system (\ref{reim1})--(\ref{reimn}) to a delay system, see, e.g., \cite{auriol}, or computing the eigenvalues of the generator associated with system (\ref{reim1})--(\ref{reimn}), see, e.g., \cite{coron book}, \cite{neves}. We proceed using the latter approach. Defining $\tilde{z}(x,t)=e^{\sigma t}\phi(x)$ and $\tilde{v}(x,t)=e^{\sigma t}\psi(x)$, $\sigma\in\mathbb{C}$, one can conclude from (\ref{reim1})--(\ref{reimn}) that the following boundary-value problem should be satisfied
\setlength{\arraycolsep}{0pt}\begin{eqnarray}
\phi'(x)&=&-\frac{\sigma}{\bar{v}}\phi(x);\quad\psi'(x)\!\!=\!\!\frac{\sigma}{c_4}\psi(x)+\frac{c_1}{c_4}e^{-\frac{c_2x}{\bar{v}}}\phi(x)\label{phi1}\\
\phi(0)&=&-L\frac{\bar{\rho}^2}{\bar{v}}\psi(0);\quad\phi\left(D\right)\!\!=\!\!-\frac{\sigma}{c_1}e^{\frac{c_2D}{\bar{v}}}\psi\left(D\right)\label{phi3}.
%\psi'(x)&=&\frac{\sigma}{c_4}\psi(x)+\frac{c_1}{c_4}e^{-\frac{c_2x}{\bar{v}}}\phi(x)\label{phi2}\\
%\phi\left(D\right)&=&-\frac{\sigma}{c_1}e^{\frac{c_2D}{\bar{v}}}\psi\left(D\right),\label{phi4}
\end{eqnarray}\setlength{\arraycolsep}{5pt}Solving (\ref{phi1}) and using the boundary conditions (\ref{phi3}), one can conclude that, in order for nontrivial solutions to (\ref{phi1}), (\ref{phi3}) to exist, $\sigma$ should satisfy the following relation
\begin{eqnarray}
-a_1e^{-\sigma\tau D}+\sigma\tau_{\rm mix}c_1-\sigma c_1\frac{1-e^{-\sigma\tau D}e^{-\frac{c_2D}{\bar{v}}}}{\sigma\tau\bar{v}+c_2}=0,\label{char0}
\end{eqnarray}
where $a_1=\frac{1}{\bar{v}}c_4c_1e^{-\frac{c_2D}{\bar{v}}}$ and $\tau=\frac{1}{c_4}+\frac{1}{\bar{v}}$. We next prove that (\ref{char0}) always admits a real solution  within an interval $\left(\sigma_1,\sigma_2\right)$ for some $\sigma_1,\sigma_2>0$. Since $c_2>0$, for $\sigma\in[0,+\infty)$ relation (\ref{char0}) is equivalent to $a_2 \sigma^2- a_1 \left(\sigma+c_2\right)e^{-\sigma\tau D}=0$, where $a_2=\bar{v}c_1\tau_{\rm mix}\tau$. The proof is completed observing that $f\left(\sigma\right)=a_2 \sigma^2- a_1 \left(\sigma+c_2\right)e^{-\sigma\tau D}$ is continuous for all $\sigma\geq0$ as well as that $f(0)=-a_1c_2<0$ and $\lim_{\sigma+\infty}f\left(\sigma\right)=+\infty$. %Thus, for any values of the parameters $a_1$, $a_2$, and $\tau$ there exists always a strictly positive real root of (\ref{char}).

\setcounter{equation}{0}
\renewcommand{\theequation}{B.\arabic{equation}}
\appendices
\section*{Appendix B: Proof of Theorem \ref{theorem1}}
%\subsection*{Proof of Theorem \ref{theorem1}}

We start defining the following functionals 
\begin{eqnarray}
V_{1_p}(t)&=&\int_0^De^{-2k_1px}\tilde{z}(x,t)^{2p}dx\nonumber\\
&&+\int_0^De^{-2k_1px}\tilde{z}_x(x,t)^{2p}dx\label{lyap1}\\
V_{2_p}(t)&=&\int_0^De^{2k_2px}\tilde{v}(x,t)^{2p}dx\nonumber\\
&&+\int_0^De^{2k_2px}\tilde{v}_x(x,t)^{2p}dx\\
V_{3_p}(t)&=&\tilde{v}\left(D,t\right)^{2p},\label{lyap3}
\end{eqnarray}
where $k_1$, $k_2$ are arbitrary positive constants and $p$ is a positive integer. Moreover, from system (\ref{tilz})--(\ref{speed2}) we obtain
\setlength{\arraycolsep}{0pt}\begin{eqnarray}
\tilde{z}_{xt}(x,t)+\bar{v}\tilde{z}_{xx}(x,t)&=&c_2\tilde{z}_x(x,t)-ke^{\frac{c_2x}{\bar{v}}}\bar{h}_{\rm mix}\bar{\rho}^2\tilde{v}_x(x,t)\nonumber\\
&&-\frac{c_2}{\bar{v}}ke^{\frac{c_2x}{\bar{v}}}\bar{h}_{\rm mix}\bar{\rho}^2\tilde{v}(x,t)\label{tilzder}\\
\tilde{v}_{xt}(x,t)-c_4\tilde{v}_{xx}(x,t)&=&-k\tilde{v}_x(x,t)\label{speed1der}\\
\tilde{z}_{x}(0,t)&=&L\frac{\bar{\rho}^2}{\bar{v}^2}c_4\tilde{v}_x(0,t)-\frac{k\bar{\rho}^2}{\bar{v}}\nonumber\\
&&\times\left(\bar{h}_{\rm mix}+\frac{L}{\bar{v}}+L\frac{c_2}{\bar{v}k}\right)\tilde{v}(0,t)\\
\tilde{v}_x(D,t)&=&0.\label{speed2der}
\end{eqnarray}\setlength{\arraycolsep}{5pt}Differentiating (\ref{lyap1})--(\ref{lyap3}) along (\ref{tilz})--(\ref{speed2}), (\ref{tilzder})--(\ref{speed2der}), with integration by parts, Young's inequality and as $\left(d_1+d_2\right)^p\leq 2^p\left(d_1^p+d_2^p\right)$, $\forall$ $d_1,d_2, p>0$, with $p$ integer, we get
\begin{eqnarray}
\dot{V}_{1_p}(t)&\leq&-\bar{v}e^{-2k_1pD}\tilde{z}\left(D,t\right)^{2p}+\left(c_7^{2p}+c_9^{2p}\right)\bar{v}\tilde{v}\left(0,t\right)^{2p}\nonumber\\
&&-2p\left(\bar{v}k_1-c_2-\frac{2p-1}{2p}c_6\right)V_{1_p}(t)+c_6V_{2_p}(t)\nonumber\\
&&-\bar{v}e^{-2k_1pD}\tilde{z}_x\left(D,t\right)^{2p}+c_8^{2p}\bar{v}\tilde{v}_x\left(0,t\right)^{2p}\label{lyap1n}\\
\dot{V}_{2_p}(t)&=&c_4e^{2k_2pD}V_{3_p}(t)-c_4\tilde{v}(0,t)^{2p}-c_4\tilde{v}_x(0,t)^{2p}\nonumber\\
&&-2p\left(k+c_4k_2\right)V_{2_p}(t)\\
\dot{V}_{3_p}(t)&=&-2pkV_{3_p}(t),\label{lyap3n}
\end{eqnarray}
with $c_6=ke^{\frac{c_2D}{\bar{v}}}\bar{h}_{\rm mix}\bar{\rho}^2\left(1+\frac{c_2}{\bar{v}}\right)$, $c_7=L\frac{\bar{\rho}^2}{\bar{v}}$, $c_8=2L\frac{\bar{\rho}^2}{\bar{v}^2}c_4$, $c_9\!=\!2\frac{k\bar{\rho}^2}{\bar{v}}\left(\bar{h}_{\rm mix}\!+\!\frac{L}{\bar{v}}\!+\!L\frac{c_2}{\bar{v}k}\right)$. Defining the Lyapunov functional
\begin{eqnarray}
V_p(t)=V_{1_p}(t)+{k_3}^{2p}V_{2_p}(t)+{k_4}^{2p}e^{2k_2pD}V_{3_p}(t),\label{deflyap}
\end{eqnarray}
we obtain from (\ref{lyap1n})--(\ref{lyap3n})
\begin{eqnarray}
\dot{V}_p(t)&\leq& -\left({k_3}^{2p}c_4-c_7^{2p}\bar{v}-c_9^{2p}\bar{v}\right)\tilde{v}\left(0,t\right)^{2p}\nonumber\\
&&-2p\left(\bar{v}k_1-c_2-\frac{2p-1}{2p}c_6\right)V_{1_p}(t)\nonumber\\
&&-\left(2pk_3^{2p}\left(k+c_4k_2\right)-c_6\right)V_{2_p}(t)\nonumber\\
&&-\left(2pk{k_4}^{2p}-{k_3}^{2p}c_4\right)e^{2k_2pD}V_{3_p}(t)\nonumber\\
&&-\left({k_3}^{2p}c_4-c_8^{2p}\bar{v}\right)\tilde{v}_x\left(0,t\right)^{2p}.\label{lyapan1}
\end{eqnarray}
Since $p\geq 1$, choosing $k_1=\frac{1}{\bar{v}}\left(c_2+c_6+k\right)$, $k_2=\frac{c_6}{2c_4}$, ${k_3}=\max\left\{ \left(c_7+c_8+c_9\right)\max\left\{\frac{\bar{v}}{c_4},1\right\},1\right\}$, and ${k_4}={k_3}\max\left\{\frac{c_4}{k},1\right\}$, we get from (\ref{lyapan1}) that $\dot{V}_p(t)\leq-2pkV_{1_p}(t)-2pk{k_3}^{2p}V_{2_p}(t)-pk{k_4}^{2p}e^{2k_2pD}V_{3_p}(t)$. Hence, using (\ref{deflyap}) we arrive at\footnote{To derive (\ref{alm}) we also used (\ref{tilzder}), (\ref{speed1der}), which implies that, in principle, higher regularity of solutions is needed. Yet, one could show that (\ref{alm}) still holds (in the sense of distributions), using similar arguments to, e.g., \cite{coron book}.} 
\begin{eqnarray}
\dot{V}_p(t)&\leq&-pk{V}_p(t).\label{alm}
\end{eqnarray}
From (\ref{alm}) we then get ${V}_p^{\frac{1}{2p}}(t)\leq e^{-\frac{k}{2} t}{V}_p^{\frac{1}{2p}}(0)$, and thus, using (\ref{deflyap}) we obtain $V_{1_p}^{\frac{1}{2p}}(t)+k_3V_{2_p}^{\frac{1}{2p}}(t)+k_4e^{k_2D}V_{3_p}^{\frac{1}{2p}}(t)\leq 4e^{-\frac{k}{2} t}\left(V_{1_p}^{\frac{1}{2p}}(0)+k_3V_{2_p}^{\frac{1}{2p}}(0)+k_4e^{k_2D}V_{3_p}^{\frac{1}{2p}}(0)\right)$. With definitions (\ref{lyap1})--(\ref{lyap3}) and standard inequalities we get
\begin{eqnarray}
\Xi_p(t)&\leq& 8e^{-\frac{k}{2} t}\Xi_p(0)\label{sh1}\\
\Xi_p(t)&=&\|\tilde{z}(t)\|_{-k_1,2p}+k_3\|\tilde{v}(t)\|_{k_2,2p}+\|\tilde{z}_x(t)\|_{-k_1,2p}\nonumber\\
&&+k_3\|\tilde{v}_x(t)\|_{k_2,2p}+k_4e^{k_2D}\left|\tilde{v}(D,t)\right|.\label{sh2}
\end{eqnarray}
Taking the limit of (\ref{sh1}) as $p\to+\infty$, with the definition of the supremum norm and (\ref{sh2}) we get $\|\tilde{z}(t)\|_{C^1}+\|\tilde{v}(t)\|_{C^1}+\left|\tilde{v}(D,t)\right|\leq\bar{\mu}\left(\|\tilde{z}(0)\|_{C^1}+\|\tilde{v}(0)\|_{C^1}+\left|\tilde{v}(D,0)\right|\right)e^{-\frac{ k}{2}t}$, for some positive constant $\bar{\mu}$, where we also used the facts that $e^{-k_1D}\|\tilde{z}(t)\|_{C}\leq    \|\tilde{z}(t)\|_{-k_1,C}\leq \|\tilde{z}(t)\|_{C}$ and $\|\tilde{v}(t)\|_{C}\leq   \|\tilde{v}(t)\|_{k_2,C}\leq e^{k_2D}\|\tilde{v}(t)\|_{C}$.  The proof is completed using relation $\tilde{z}(x)=e^{\frac{c_2x}{\bar{v}}}\left(\tilde{\rho}(x)+\bar{h}_{\rm mix}\bar{\rho}^2\tilde{v}(x)\right)$ .

%The fact $\tilde{z}(x)=e^{\frac{c_2x}{\bar{v}}}\left(\tilde{\rho}(x)+\bar{h}_{\rm mix}\bar{\rho}^2\tilde{v}(x)\right)$ completes the proof.

% The proof is completed using relation $\tilde{z}(x)=e^{\frac{c_2x}{\bar{v}}}\left(\tilde{\rho}(x)+\bar{h}_{\rm mix}\bar{\rho}^2\tilde{v}(x)\right)$ .

\setcounter{equation}{0}
\renewcommand{\theequation}{C.\arabic{equation}}
\appendices
\section*{Appendix C: Proof of Theorem \ref{theorem2}}
%\subsection*{Proof of Theorem \ref{theorem2}}

%First note that estimate (\ref{thm1est}) implies the existence of the temporal norms incorporated in inequalities (\ref{est st}), (\ref{est st1}) for all $x_1,x_2\in[0,D]$. 
Taking the Laplace transform of (\ref{speed1}) and setting the initial condition to zero we obtain
\begin{eqnarray}
\tilde{V}\left(x,s\right)=e^{-\frac{k}{c_4}\left(x_1-x\right)}e^{-\frac{s}{c_4}\left(x_1-x\right)}\tilde{V}\left(x_1,s\right),\label{diffv}
\end{eqnarray}
for all $0\leq x\leq x_1\leq D$, and thus, 
\begin{eqnarray}
\tilde{v}\left(x_2,t\right)=e^{-\frac{k}{c_4}\left(x_1-x_2\right)}\tilde{v}\left(x_1,t+\frac{x_2-x_1}{c_4}\right),\label{p1}
\end{eqnarray}
where $\tilde{v}\left(x_1,\theta\right)$, within the interval $\frac{x_2-x_1}{c_4}\leq \theta<0$, is set to zero, as, in the present context, we are concerned with an input-output representation (without accounting for the effect of initial conditions), in which $\tilde{v}\left(x_1,t\right)$ is viewed as input signal (with initial condition $\tilde{v}\left(x_1,\theta\right)$, $\frac{x_2-x_1}{c_4}\leq \theta<0$). % (as in the present, input-output framework, $\tilde{v}\left(x_1,t\right)$ is viewed as input signal with zero initial condition $\tilde{v}\left(x_1,\theta\right)$, $\frac{x_2-x_1}{c_4}\leq \theta<0$). 
Therefore, from (\ref{p1}) we get that $\left(\int_0^{+\infty}\left|\tilde{v}\left(x_2,t\right)\right|^{p}dt\right)^{\frac{1}{p}}=e^{-\frac{k}{c_4}\left(x_1-x_2\right)}\left(\int_0^{+\infty}\left|\tilde{v}\left(x_1,t+\frac{x_2-x_1}{c_4}\right)\right|^{p}dt\right)^{\frac{1}{p}}$. Since $k,c_4>0$ and $x_2<x_1$, we obtain (\ref{est st}) for $p\in[1,+\infty)$. Similarly, taking a supremum over time in (\ref{p1}) we obtain (\ref{est st}) for $p=+\infty$. Differentiating (\ref{diffv}) with respect to $x$ we get $\tilde{V}'\left(x,s\right)=e^{-\frac{k}{c_4}\left(x_1-x\right)}e^{-\frac{s}{c_4}\left(x_1-x\right)}\tilde{V}'\left(x_1,s\right)$, for all $0\leq x\leq x_1\leq D$. Thus, employing identical arguments we obtain (\ref{est st1}).

\section*{Acknowledgments}
Nikolaos Bekiaris-Liberis was supported by the funding from the European Commission's Horizon 2020 research and innovation programme under the Marie Sklodowska-Curie grant agreement No. 747898, project PADECOT. 

The authors thank Dr. Diamantis Manolis, Prof. Markos Papageorgiou, and Prof. Claudio Roncoli for fruitful discussions.

%\section{Adaptive Control Design}

%\section{Trajectory Generation for the Linearized System}

%\section{In-Domain Disturbance Rejection Not So Good}
%We introduce the following integral state
%\begin{eqnarray}
%\sigma_t(x,t)=\frac{1}{\tau_{\rm mix}}\left(\frac{1}{h_{\rm mix}\left(h_{\rm acc}(x,t)\right)}\left(\frac{1}{\rho(x,t)}-L\right)-v(x,t)\right).
%\end{eqnarray}
%Its linearization is given by
%\begin{eqnarray}
%\tilde{\sigma}_t(x,t)=\frac{1}{\tau_{\rm mix}}\left(-\frac{1}{\bar{h}_{\rm mix}\bar{\rho}^2}\tilde{\rho}(x,t)-\tilde{v}(x,t)-\frac{h'\left(\bar{h}_{\rm acc}\right)}{\bar{h}_{\rm mix}^2}\left(\frac{1}{\bar{\rho}}-L\right)\tilde{h}(x,t)\right),
%\end{eqnarray}
%and hence,
%\begin{eqnarray}
%\tilde{\sigma}_t(x,t)=-c_1e^{-\frac{c_2x}{\bar{v}}}\tilde{z}(x,t)-c_3\tilde{h}(x,t).
%\end{eqnarray}
%So in closed-loop 
%\begin{eqnarray}
%\tilde{\sigma}_t(x,t)=-k\tilde{v}(x,t)-k_1\tilde{\sigma}(x,t).
%\end{eqnarray}
%delay in $h$ and trajerctory trcking ?

%why steady-state profile depends on relaxation time as well? Reasonable if the time gaps, or the equilibria for speed, are different. 

%why this FD? For cth strategy 

%\markboth{Submitted to Transportation Research Part B on March 6, 2017}{Bekiaris-Liberis, Roncoli, and Papageorgiou}

%\begin{abstract}
%\baselineskip=1.7\normalbaselineskip

%\end{abstract}
%\baselineskip=1.8\normalbaselineskip

\end{document}